\documentclass[11]{amsart}

\usepackage{amsfonts}
\usepackage{amscd}
\usepackage{amsmath, mathrsfs, amssymb}
\usepackage{amsthm}
\usepackage{setspace}
\usepackage{hyperref}
\usepackage{color}
\usepackage{epsfig}
\usepackage{here}
\usepackage{graphicx}
\usepackage[all]{xy}
\usepackage{psfrag}
\usepackage{graphicx,transparent}
\usepackage{enumerate}
\usepackage{caption}

\theoremstyle{plain}

\theoremstyle{definition}

\newcommand{\MM}{\mathcal M}

\newcommand{\PP}{\mathcal P}

\newcommand{\OO}{\mathcal O}

\newcommand{\HH}{\mathcal H}

\newcommand{\bbC}{\mathbb C}

\begin{document}

\title[]{Complete curves in the strata of differentials}

\date{\today}

\author{Dawei Chen}
      \thanks{Research partially supported by National Science Foundation Grant DMS-2001040
        and Simons Foundation Collaboration Grant 635235.}
\address{Department of Mathematics, Boston College, Chestnut Hill, MA 02467, USA}
\email{dawei.chen@bc.edu}

\begin{abstract}
Gendron proved that the strata of holomorphic differentials with prescribed orders of zeros do not contain complete algebraic curves by applying the maximum modulus principle to saddle connections. Here we provide an alternative proof for this result by using positivity of divisor classes on moduli spaces of curves.
\end{abstract}

\maketitle


For $\mu = (m_1, \ldots, m_n)$ with $\sum_{i=1}^n m_i  = 2g-2$ and $m_i \geq 0$ for all $i$, let $\HH_g(\mu)$ be the stratum of holomorphic differentials $\omega$ of zero type $\mu$ on smooth and connected complex genus $g$ curves $X$, i.e., the underlying canonical divisor of $\omega$ is of type $\sum_{i=1}^n m_i p_i$ for distinct points $p_1, \ldots, p_n$ in $X$.

The study of differentials is significant in surface dynamics and moduli theory.  We simply refer to \cite{Zorich, Wright, ChenSurvey} for an introduction to this fascinating subject. An important question for understanding the geometry of a (non-compact) moduli space is to study complete subvarieties contained in it. For instance for the moduli space of curves $\MM_g$, this amounts to finding families of smooth curves that have non-constant complex structures and do not degenerate when approaching to the boundary. While there exist complete subvarieties of positive dimension in $\MM_g$, the sharp upper bound for their dimensions is largely unknown (see \cite{Diaz} for an upper bound), e.g., it is unknown whether $\MM_4$ contains a complete (algebraic) surface. 

One can similarly ask this question for the strata of holomorphic differentials. In \cite{Gendron} Gendron proved the following nice result. 

\subsection*{Theorem}{\em 
For $\mu = (m_1, \ldots, m_n)$ with $m_i \geq 0$ for all $i$, the stratum $\HH_g(\mu)$ of holomorphic differentials of type $\mu$ does not contain complete curves. } 

\bigskip

Gendron's argument relies on flat geometry. It applies the maximum modulus principle to saddle connections (i.e., geodesics joining the zeros of $\omega$ under the induced metric), which essentially forces a complete family of differentials to be trivial. Below we provide a short algebraic proof for this result, which uses positivity of divisor classes on moduli spaces of curves. 

We first review some related moduli spaces and divisor classes. Denote by $\PP_g(\mu) = \HH_g(\mu)/\bbC^{*}$ the projectivized stratum parameterizing canonical divisors of type $\mu$; then $\HH_g(\mu)$ can be identified with the tautological bundle $\OO(-1)$ over $\PP_g(\mu)$ with the zero section removed. We denote by $\eta = c_1 (\OO(-1))$ the first Chern class of the tautological bundle. Consider the zeros of the canonical divisors of type $\mu$ as marked points in the underlying curves; then $\PP_g(\mu)$ can be embedded in $\MM_{g,n}$ as a subvariety. Let $\kappa = f_{*}(c_1(\omega_{f})^2)$ where $\omega_f$ is the relative dualizing sheaf on the universal curve $f\colon \mathcal X\to \MM_{g,n}$. 
Let $\psi_i$ be the first Chern class of the bundle of cotangent lines with respect to the $i$-th marked point. We also use the same notations for the restrictions of these divisor classes from $\MM_{g,n}$ to $\PP_g(\mu)$. It is known that the divisor classes $\psi_i$ and $\kappa$ are proportional to $\eta$ on $\PP_g(\mu)$ (see \cite[Proposition 2.1]{ChenTauto}): 
$$ \eta = (m_i+1)\psi_i, \quad \kappa = \kappa_\mu \eta $$ 
where $\kappa_\mu = 2g-2+n - \sum_{i=1}^n 1 / (m_i+1)$.  

\begin{proof}[Proof of the Theorem]
Suppose $C$ is a complete (irreducible) curve contained in $\HH_g(\mu)$.   Let $C_0\subset \PP_g(\mu)$ be the image of $C$ under the projection $\HH_g(\mu)\to \PP_g(\mu)$. Since each fiber of $\HH_g(\mu)$ over $\PP_g(\mu)$ is $\bbC^{*}$, the complete curve $C$ cannot be contained in a fiber.  Therefore, $C$ corresponds to a multi-section of $\OO(-1)$ over $C_0$. Let $s$ be the degree of $C$ over $C_0$. Then $C$ gives a nowhere vanishing section of $\OO(-s)$ over $C_0$, which implies that $s\eta$ and hence $\psi_i$ and $\kappa$ are all trivial restricted to $C_0$ by the above divisor class relations. However, $\kappa + \sum_{i=1}^n \psi_i$ is an ample divisor class on $\overline{\MM}_{g,n}$ (see \cite[Theorem (2.2)]{Cornalba} and \cite[Chapter XIV, Theorem (5.1)]{ACG})\footnote{This ample divisor class can also be written as $12\lambda - \delta + \sum_{i=1}^n \psi_i$ by using the well-known relation $12\lambda = \kappa + \delta$ where $\lambda$ is the first Chern class of the Hodge bundle and $\delta$ is the total boundary divisor class. Note that our $\kappa$ class differs from the one in \cite{ACG} by the total $\psi$ class, as we use the ordinary relative dualizing sheaf instead of the log version in the definition of $\kappa$.} which thus has positive degree on every complete curve contained in it and cannot be trivial on $C_0$. This leads to the desired contradiction. 
\end{proof}

\subsection*{Remark}
For the strata of meromorphic differentials, i.e., when $\mu$ contains some negative entries, in \cite{ChenAffine} it was shown that both $\HH_g(\mu)$ and $\PP_g(\mu)$ do not contain complete curves. The proof relies on the sign change of $m_i +1$ at a zero and at a pole. However for the case of all $m_i \geq 0$, a complete curve $C_0$ in $\PP_g(\mu)$, if it exists, would not necessarily lift to a complete curve in $\HH_g(\mu)$, e.g., the lift might cross the zero section of $\OO(-1)$. Hence it remains open to determine whether $\PP_g(\mu)$ contains a complete curve when all the entries of $\mu$ are nonnegative. 

Note that if one adds the boundary divisor $\Delta_0$ parameterizing differentials on nodal curves with non-separating nodes only, then for every signature $\mu$ with nonnegative entries the partial compactification  $\PP_g(\mu) \cup \Delta_0$ contains infinitely many complete curves arising from Teichm\"uller curves (see \cite[Corollary 3.2]{ChenMoeller}). 

Finally one can study the strata of $k$-differentials $\HH^k_g(\mu)$ and $k$-canonical divisors $\PP^k_g(\mu)$ where $\mu$ is a partition of $k(2g-2)$. 
  Via the canonical $k$-cover $\pi$, a $k$-differential $\xi$ can lift to an abelian differential $\omega$ (up to $k$-th roots of unity), i.e., $\pi^{*}\xi = \omega^k$, 
  where $\omega$ is holomorphic if and only if $\xi$ has no pole of order $\leq -k$ (see e.g., \cite[Section 2.1]{BCGGM-k}). Therefore, $\HH^k_g(\mu)$ contains no complete curves for any $\mu$ and $\PP^k_g(\mu)$ contains no complete curves if $\mu$ has an entry $\leq -k$.  It remains open to determine whether $\PP^k_g(\mu)$ contains a complete curve when all entries of $\mu$ are bigger than $-k$.  

\subsection*{Acknowledgements} The author thanks the referees for carefully reading the paper and very helpful comments.  
  

\end{document}